\documentclass[oribibl]{llncs}

\usepackage{amsfonts,amsmath,amssymb,amsbsy,bm}
\usepackage{amsfonts}
\usepackage{graphicx}
\usepackage{amssymb}
\usepackage{epstopdf}
\usepackage{graphicx}
\usepackage{cite}
\usepackage{hyperref}
\usepackage[usenames]{color}

\definecolor{ascolor}{rgb}{0,0.6,0}
\definecolor{mblcolor}{rgb}{0,0,0.8}
\definecolor{daocolor}{rgb}{0.4,0,0.25}
\definecolor{pbbcolor}{rgb}{0.8,0,0}

\usepackage[normalem]{ulem}

\DeclareMathAccent{\maxvec}{\mathord}{letters}{"7E}
\DeclareMathOperator*{\argmin}{arg\,min}

\usepackage{subfigure}
\usepackage{cleveref}
\spnewtheorem{conj}[theorem]{Conjecture}{\bfseries}{\itshape}

\numberwithin{equation}{section}
\newcommand{\bbR}{\mathbb{R}}  
\newcommand{\bbZ}{\mathbb{Z}}  

\newcommand{\calE}{\mathcal{E}} 
\newcommand{\calL}{\mathcal{L}} 
\newcommand{\calR}{\mathcal{R}} 
\newcommand{\calU}{\mathcal{U}} 

\newcommand{\mF}{\mathsf{F}}
\newcommand{\mG}{\mathsf{G}}

\renewcommand{\a}{{\rm a}}
\renewcommand{\c}{{\rm c}}
\renewcommand{\o}{{\rm o}}
\newcommand{\atc}{{\rm atc}}
\newcommand{\core}{{\rm core}}

\newcommand{\Ea}{\calE^\a}

\newcommand{\Ua}{\calU^\a}
\newcommand{\Uaz}{\calU^\a_0}
\newcommand{\Ec}{\calE^\c}

\newcommand{\Uc}{\calU^\c}
\newcommand{\Ucz}{\calU^\c_0}

\begin{document}

\frontmatter
\title{Development of an Optimization-Based Atomistic-to-Continuum Coupling Method}

\author{Derek Olson\inst{1} \and Pavel Bochev\inst{2} \and Mitchell Luskin\inst{1} \and Alexander V. Shapeev\inst{1}}

\institute{University of Minnesota, Minneapolis, United States \\ \email{\{olso4056, luskin, ashapeev\}@umn.edu}\thanks{DO was supported by the Department of Defense (DoD) through the National Defense Science \& Engineering
Graduate Fellowship (NDSEG) Program.  ML was supported in part by the NSF PIRE Grant OISE-0967140,
DOE Award DE-SC0002085, and AFOSR Award
FA9550-12-1-0187. AS was supported in part by the
DOE Award DE-SC0002085 and AFOSR Award
FA9550-12-1-0187.} \and Sandia National Laboratories, Albuquerque, United States \\ \email{pbboche@sandia.gov}\thanks{
		 Sandia National Laboratories is a multi-program laboratory
                managed and operated by Sandia Corporation, a wholly owned subsidiary of
                Lockheed Martin Corporation, for the U.S. Department of
                Energy's National Nuclear Security Administration under
                contract DE-AC04-94AL85000.}}

\maketitle

\begin{abstract}
Atomistic-to-Continuum (AtC) coupling methods are a novel means
of computing the properties of a discrete crystal structure,
such as those containing defects, that combine the accuracy of
an atomistic (fully discrete) model with the efficiency of a
continuum model. In this note we extend the
optimization-based AtC, formulated in~\cite{olson} for linear,
one-dimensional problems to multi-dimensional settings and
arbitrary interatomic potentials.
We conjecture optimal error estimates for the multidimensional AtC, outline an implementation procedure, and provide numerical results to corroborate the conjecture for a 1D Lennard-Jones system with next-nearest neighbor interactions.
\end{abstract}

\pagestyle{headings}

\section{Introduction}\label{S:intro}

Solid materials have atomic configurations which are arranged as a crystalline lattice, and the properties of these materials are derived from the underlying structure of the lattice.  Specifically, defects in the regular, repeating arrangement of atoms such as a dislocation, or an extra plane of atoms, determine fundamental mechanisms such as plastic slip.  The presence of defects invalidate the central hypotheses of continuum mechanics so models that recognize the discrete nature of the material on the atomic scale must be used.  Such methods can vary in their complexity ranging from quantum mechanical models which incorporate nuclear and electronic forces to empirical potential models that assume the existence of a potential energy which is a function of the nuclear positions only.  The latter allows atoms to be considered as classical mechanical particles.  Throughout this note, we assume that the exact mathematical problem we wish to solve is that of minimizing the global potential energy of a set of $N$ atoms or, equivalently, of equilibrating the internal and external forces on the atoms.

The outstanding issue with empirical atomistic models is the complexity involved in their applications.  In even the smallest problems of material interest on the nanoscale, there will be at least $10^9$ and up to $10^{15}$ atoms meaning the number of degrees of freedom in an atomistic model is often far outside the scope of any current computational feasibility.  A novel attempt at solving this problem has been to keep the atomistic model only in a small region near the defect, while employing a continuum model such as elasticity in the bulk of the material away from the defect. Continuum models are well understood and can numerically be solved in an efficient manner using finite elements.  In effect, the atomistic model provides a constitutive relation near the defect where the constitutive relation of the continuum model fails to hold.

These so called atomistic-to-continuum (AtC) coupling methods have seen a surge of interest in the last two decades, especially with the introduction of the quasicontinuum method in~\cite{tadmor1996}.  The problem introduced in these AtC methods is how to combine, or \textit{couple}, the two different models.  An informal way of carrying this out is to divide the computational domain, say $\Omega$, into an atomistic region, $\Omega_\a$, and a continuum region, $\Omega_\c$.  Then, a global hybrid energy or hybrid force field is constructed from the atomistic and continuum models on $\Omega_\a$ and $\Omega_\c$.  The resulting hybrid energy is then minimized, or alternatively, the internal and external forces are equilibrated to find the equilibrium configuration of $\Omega$.

In this note we continue the development of the optimization-based AtC approach commenced in~\cite{olson}.  The core idea is to pose  independent atomistic and continuum subproblems on overlapping domains $\Omega_\a$ and $\Omega_\c$ and then couple the models by minimizing an objective functional, which measures the difference between the strains of the atomistic and continuum states on $\Omega_\a \cap \Omega_\c$.
In so doing, our approach combines ideas from \emph{blending AtC} methods~\cite{bqcf,Bochev_08_MMS,belytschko2003coupling,koten_2011,Bauman_08_CM,bqcf13,bqce12} with the optimization-based domain-decomposition approach for PDEs in \cite{Gunzburger_00_CMA,Gunzburger_00_SINUM}.

The resulting optimization-based AtC method differs substantially from current energy or force-based methods, and to the best of our knowledge \cite{olson} is the first instance of using an objective functional of this form to effect atomistic-to-continuum coupling. Conceptually, our AtC approach is similar to the heterogeneous domain decomposition method for PDEs developed in \cite{Gervasio_01_NM} with the important distinction that we couple two fundamentally different material models rather than PDEs.

The main focus of this note is on the formulation of an
optimization-based AtC method for modeling material defects in two and three dimensions, while allowing for
arbitrary many-body terms in the potential energy.
Sect.~\ref{S:prelim} quotes the necessary background
results and Sect.~\ref{S:formulate} presents the formulation
of the method. Solution of the optimization problem is
discussed in Sect.~\ref{S:imp}. We conjecture error
estimates and derive optimal parameters for our algorithm from
the complexity analysis of Sect.~\ref{S:error}. Finally,
Sect.~\ref{S:num} provides numerical evidence in support of
these conjectures.

\section{Preliminaries}\label{S:prelim}

\begin{figure}[htp!]
\centering
{\resizebox{2in}{2in}{\includegraphics{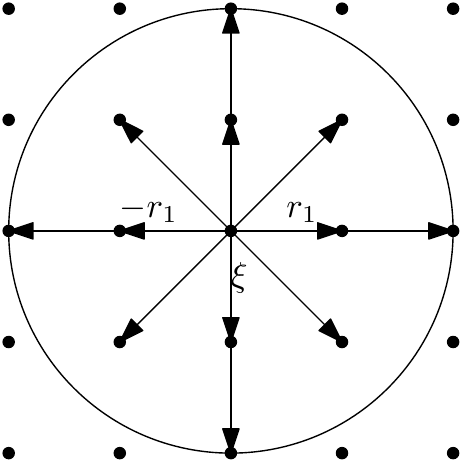}}}
\caption{An atom site $\xi$ and its interaction range $\mathcal{R}$.}
\label{fig:manyInt}
\end{figure}

We consider the problem of modeling a crystal occupying the infinite domain, $\mathbb{R}^d$, and take the reference configuration of the atoms to be the integer lattice, $\bbZ^d$, deformed by the macroscopic deformation gradient $\mF$.  Deformations of the material are thus described by functions $y:\mF\mathbb{Z}^d \to \mathbb{R}^d$.  For any deformed configuration, $y$, of the lattice, we assume the energy due to electronic and nuclear interactions can be described by an empirical site potential $V_\xi(y)$ where $V_\xi$ represents the energy attributable to atom $\xi \in \mathbb{Z}^d$.  As usual, we further assume that each $\xi$ interacts with only a finite number of other atoms.  The set of atoms that $\xi$ interacts with is given by $\xi + \calR \subset \bbZ^d$ where $\calR$ is the interaction neighborhood.  The interaction neighborhood can be defined through a cutoff radius, $r_{\rm cut}$, so that
\[
\mathcal{R} =~ \left\{ \xi \in \mathbb{Z}^d \, | \, 0 < |\mF \xi| \leq r_{\rm cut} \right\}.
\]
Fig.~\ref{fig:manyInt} depicts $\mathcal{R}$ in $2D$ where $\mF$ is the identity and $r_{\rm cut} = 2$.
We model point defects in the lattice by allowing $V_\xi$ to
depend on $\xi$ while assuming that $V_\xi = V$ when $\xi$ is
far from a defect. An evident example is an impurity where
atoms of a different species have different interaction laws
with the bulk atoms, but these interactions are only limited to
small neighborhoods of defect (impure) atoms.

The presence of defects in the lattice generates elastic
fields causing the atoms to relax. The deformed configurations
are generically given by
\[
y(\xi) = \mF\xi + u(\xi),
\]
where $u:\bbZ^d \to \bbR^d$ is the displacement field.  The
energy of the deformed configuration associated to this
displacement field is
\begin{equation}\label{enDef}
\calE(u) := \sum_{\xi \in \bbZ^d}
	V_\xi(Du(\xi))
\end{equation}
where $Du(\xi) := (D_\rho u(\xi))_{\rho\in\calR}$ is a collection of finite differences of $u$ and $D_{\rho}
u(\xi) := u(\xi + \rho) - u(\xi)$ defines the finite difference operator. We note that $V_\xi$
implicitly depends on the macroscopic deformation gradient
$\mF$. Furthermore, without loss of generality, we assume that
$V_\xi(0)=0$ so that the infinite sum in $\calE(u)$ is
well-defined. For example, in the case of a Lennard-Jones
potential, $\phi$, with next-nearest neighbor interactions in
$1D$, we can define
\[
V_\xi(Du(\xi)) = \phi(\mF + D_1u) + \phi(2\mF + D_1u - D_{-1}u) - (\phi(\mF) + \phi(2\mF))
\]
where $\phi(\mF) + \phi(2\mF)$ is subtracted from the
usual Lennard-Jones potential (without affecting the computed
forces) so that $V_\xi(0)=0.$

 The problem we seek to solve is then
\begin{equation}\label{atProb}
\bar{u} \in \argmin_{u\in\mathcal{U}} \calE(u),
\end{equation}
where $\argmin$ denotes the set of local minima of a functional and the admissible displacement space is taken to be $\mathcal{U} = \left\{u: \bbZ^d \to \bbR^d\right\}$.  Typically, this energy on an infinite domain is approximated by truncating to a finite domain (the approach taken here) or by imposing periodic boundary conditions.  However, the complexity involved in computing the resulting energy may be intractable for current computing capabilities due to the large number of atoms and interactions so a more efficient stratagem is required.

One solution approach would be to use continuum hyperelasticity models, but the elastic fields involved in modeling defects such as dislocations are singular at the defect core and so do not belong to the function spaces required in  a standard continuum formulation.  Atomistic-to-continuum models seek to overcome these deficiencies by utilizing both models simultaneously: the atomistic model near the defect and the continuum model far from the defect.

\section{An AtC Method Formulation} \label{S:formulate}

\subsection{Decomposition into atomistic and continuum subdomains}

Typical AtC methods require the decomposition of the computational domain $\Omega$ into atomistic and continuum subdomains, $\Omega_\a$ and $\Omega_\c$, respectively, with a possible blending, or overlap, region $\Omega_\o := \Omega_\a \cap \Omega_\c$.  The goal of these methods is to create a \textit{globally} defined hybrid energy or force field derived from using the atomistic model in $\Omega_\a$, the continuum model in $\Omega_\c$, and some coupling of the two in $\Omega_\o$.  The distinguishing feature of our algorithm is to pose the atomistic and continuum problems independently on overlapping domains and then couple them by minimizing a suitably defined norm of the difference between the separate atomistic and continuum states that exist simultaneously on the overlap region.  As we shall see, some care must be taken in the definitions of $\Omega_\a$ and $\Omega_\c$ to account for the interaction range, $\mathcal{R}$, from the previous section.

Truncation of the infinite domain, $\mathbb{R}^d$, to a finite, regular polygonal domain, $\Omega$, inscribed in a sphere of radius $R_\c$, is the first approximation in  modeling~\eqref{atProb}.
The boundary of $\Omega$ coincides with $\mF\mathbb{Z}^d$, and the lattice corresponding to $\Omega$ is defined as $\mathcal{L} := \Omega \cap \mF\mathbb{Z}^d$.  Consequently, we denote the space of admissible displacements which satisfy the far-field boundary condition $u(\xi) = 0$ whenever $\xi \notin \mathcal{L}$ by
\[
\mathcal{U}_0 := \left\{ u \in \mathcal{U} \, | \, u(\xi) = 0 \,\, \forall \, \xi \notin \mathcal{L}\right\}
\]
and replace~\eqref{atProb} by
\begin{equation}\label{trunProb}
\bar{u} \in \argmin_{u\in\mathcal{U}_0} \calE(u).
\end{equation}

\begin{remark}\label{rem:finiteDomain}
Though we have derived~\eqref{trunProb} with the idea of approximating an infinite domain containing a defect, a second problem of practical interest is minimizing an energy $\mathcal{E}$ on a fixed domain, $\Omega$, subject to some prescribed boundary conditions on $\partial \Omega$ and an imposed external force in $\Omega$.  In this case, we typically separate $V_\xi$ into an internal site energy $V_\xi^{\rm int}$ and an external site energy $V_\xi^{\rm ext}$.  Aside from this notational convenience, the formulation of our AtC method is identical for both of these problems.
\end{remark}

\begin{figure}[htp!]
\centering
\includegraphics[width=0.7\textwidth]{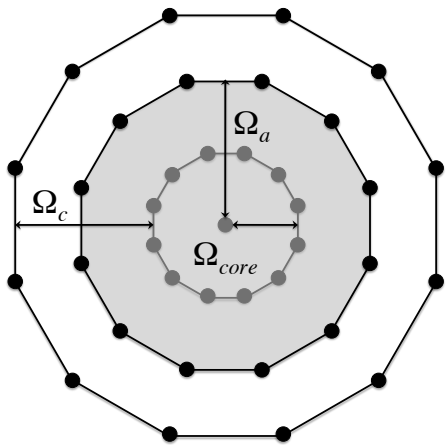}
\vspace{-2ex}
\caption{Decomposition of $\Omega$ into atomistic and continuum subdomains.}
\label{fig:domains}
\end{figure}

\begin{remark}\label{rem:domain}
For any domain, $\Omega_{\rm t} \subset \mathbb{R}^d$, (${\rm t} = \a,\c,\o$, etc.) we define its (outer) radius, $R_{\rm t} := \frac{1}{2} {\rm diam}(\Omega_{\rm t})$, and its associated discrete lattice, $\calL_{\rm t} = \Omega_{\rm t} \cap \mF\bbZ^d$.
\end{remark}
We further decompose $\Omega$ into overlapping atomistic and continuum subdomains, $\Omega_\a$ and $\Omega_\c$, as follows.  Let $\Omega_\a \subset \Omega$ be a regular polytope of radius $R_\a$ with $R_\a \ll R_\c$, and take $\Omega_\core$ to be another regular polytope of radius $R_\core < R_\a$.  The continuum subdomain is defined by $\Omega_\c := \Omega \backslash \Omega_\core^{\circ}$.  This decomposition results in an annular overlap region $\Omega_\o := \Omega_\a \cap \Omega_\c$ with width $R_\a - R_\core$.  See Fig.~\ref{fig:domains} for an illustration in 2D.
The atomistic interior of $\Omega_\a$, denoted by $\Omega_\a^{\circ}$, is the set of atoms $\xi \in \Omega_\a$ such that all neighbors of $\xi$ are also in $\Omega_\a$.
Thus
\begin{equation}\label{atCore}
\Omega_\a^{\circ} := \{\xi \in \Omega_\a \, | \,  \xi + \calR \subset \Omega_\a \}
\quad \text{and} \quad
\Omega_\a^{\circ\circ} := \{\xi \in \Omega_\a^\circ \, | \,  \xi + \calR \subset \Omega_\a^\circ \},
\end{equation}
where $\Omega_\a^{\circ\circ}$ can be interpreted as the atomistic interior of $\Omega_\a^\circ$.

For these domains, we define the associated displacement spaces
\[
\Ua := \{u:\calL_\a\to\bbR^n\}
\quad\text{and}\quad
\Uaz := \{u\in\Ua \, | \,  u = 0 \text{ outside } \Omega_\a^{\circ\circ}\}.
\]
The energy on these spaces is
\begin{equation}\label{atEnergy}
\Ea(u) :=
\sum_{\xi \in \calL_\a^{\circ\circ}}
	V_\xi(Du(\xi)),
\end{equation}
where $\calL_\a^{\circ\circ} := \Omega_\a^{\circ\circ} \cap \mathcal{L}$.  The problem of finding local minima of this energy in the space $\mathcal{U}^\a$ subject to some prescribed boundary values on $\calL_\a \backslash \calL_\a^{\circ\circ}$ is exactly what has been described as the atomistic model on $\Omega_\a$.

\begin{remark}\label{rem:impose}
As previously mentioned, some care must be exercised in defining $\Omega_\a$ and $\Omega_\c$. Precisely, we must impose the requirement that $\Omega_\core \subset \Omega_\a^{\circ\circ}$.  This ensures that the overlap width is at least twice the size of the interaction range; a necessary condition is $R_\a - R_\core \geq 2r_{\rm cut}$.
\end{remark}

Our next task is to define a continuum model on $\Omega_\c$ which is accomplished by defining the Cauchy-Born continuum energy there.  We momentarily assume a finite element triangulation, $\mathcal{T}_h$, is given in $\Omega_\c$.  This triangulation will be explicitly constructed in Sect.~\ref{S:param}.  Piecewise linear continuous finite elements are employed, and the mesh is fully refined in $\Omega_\o$ so that a finite element node exists at each $\xi \in \mathcal{L}_\o$.  We denote the space of finite elements as $\Uc$ while the subspace of $\Uc$ satisfying homogeneous Dirichlet boundary conditions on the ``outer'' boundary of $\Omega_\c$ is
\[
\Ucz := \left\{u \in \Uc \, | \, u = 0 \, \text{on} \, \partial\Omega_\c\backslash\partial\Omega_{\core}  \right\}.
\]
The Cauchy-Born continuum approximation on $\Omega_\c$ is then
\begin{equation}\label{multiCB}
\Ec(u) := \int\limits_{\Omega_\c} W(\nabla u)\, dx,
\end{equation}
where the Cauchy-Born strain energy density functional is $W(\mG) := V(\mF\calR+\mG \calR)$.  This energy is evaluated at elements of the space $\Ucz$ so that we may write the continuum energy as
\begin{equation}\label{discreteCB}
\Ec(u) = \sum_{T \in \mathcal{T}_h}  |T| \cdot  W(\nabla u|_{T}),
\end{equation}
and the continuum model consists of finding local minima in $\Uc$ of this functional subject to prescribed boundary conditions on $\partial\Omega_\c$.\footnote{In both the atomistic and continuum model, we have referenced some unknown, prescribed boundary values.  These can be interpreted as virtual controls as defined in~\cite{gervasio_2001} and discussed in~\cite{olson}.}

\subsection{Coupling}\label{S:coupling}
Having decomposed the computational domain into atomistic and continuum constituencies, we need to provide the mechanism by which these two models are coupled together.  This is done by minimizing the energy norm difference between atomistic and continuum states resulting from the atomistic and continuum problems from the spaces $\mathcal{U}^\a$ and $\mathcal{U}^\c$.  Since an atomistic state, $u^\a \in \mathcal{U}^\a$, is a discrete function defined on a lattice, whereas a continuum state, $u^\c \in \mathcal{U}^\c$, is a continuous function, we define a continuous, piecewise linear nodal interpolant of $u^\a$ on $\mathcal{T}_h$ restricted to $\Omega_\o$ by $Iu^\a$, which allows us to compare the atomistic and continuum states in the same function space on $\Omega_\o$.

Our AtC method is to then solve the constrained minimization problem
\begin{equation}\label{atcOpt}
\begin{array}{c}
\text{find $(\bar{u}^\a,\bar{u}^\c)$ such that $ \|\nabla Iu^\a - \nabla u^\c\|_{L^2(\Omega_\o)}$ is minimized} \\[1ex]
\text{subject to}
\left\{
\begin{array}{l}
\langle \delta\Ea(u^\a), v^\a\rangle=0 ~\forall v^\a\in\Uaz \\ [2ex]
\langle\delta\Ec(u^\c), v^\c\rangle=0 ~\forall v^\c\in\Ucz
\end{array}
\right.
\quad
\text{and}
\quad
\displaystyle
\int\limits_{\Omega_\o} \left(Iu^\a - u^\c\right) \, dx = 0
\end{array}
\end{equation}
The objective in (\ref{atcOpt}) ensures that the mismatch between $\bar{u}^\a$ and $\bar{u}^\c$ over $\Omega_\o$ is as small as possible. The first two constraints in (\ref{atcOpt}) imply that $\bar{u}^\a$ and $\bar{u}^\c$ are equilibria of the atomistic and continuum subproblems defined on $\Omega_\a$ and $\Omega_\c$.  The third constraint is necessary because the objective is a difference of two gradients, and without it the optimal solution would be determined only up to an arbitrary constant\footnote{In one dimension, or when there are multiple overlap regions associated with modeling multiple defects, a constraint is specified for each individual overlap region.}.
Finally, we define our AtC approximation by
\begin{equation}\label{atcApprox}
\bar{u}^\atc(x) =~ \begin{cases} \bar{u}^\a(x),\, \quad &|x| \leq R_\a, \\
\bar{u}^\c(x),\, \quad &|x| > R_\a.
\end{cases}
\end{equation}

\section{Solution of the AtC optimization problem}\label{S:imp}
The AtC formulation (\ref{atcOpt}) is a constrained
optimization problem. A standard solution approach for such
problems is to recast them into unconstrained optimization
problems through the Lagrange multiplier method. Setting the
first variations of the resulting Lagrangian with respect to
the states and the adjoints to zero yields an optimality system
from which we can determine the optimal solution of the
original problem. This approach is know as a ``one-shot
method'' \cite{Gunzburger_02_BOOK} because we solve
simultaneously for the states, adjoints, and controls.

For the AtC formulation (\ref{atcOpt}), we introduce the
Lagrange multipliers (adjoint variables) $\lambda_\a \in
\mathcal{U}^\a_{0}$ and $\lambda_\c \in \mathcal{U}^\c_{0}$
for the first two constraints, the multiplier $\eta \in
\mathbb{R}$ for the third constraint, and the Lagrangian
functional
\begin{equation}\label{lagrangian1}
\begin{array}{l}
\displaystyle
\Psi(u^\a, u^\c, \lambda_\a, \lambda_\c, \eta) = \frac{1}{2}\|\nabla Iu^\a - \nabla u^\c\|^{2}_{L^2(\Omega_\o)}  \\[1ex]
\displaystyle
\qquad\qquad
- \left< \delta \mathcal{E}^\a(u^\a), \lambda_\a\right> - \left< \delta \mathcal{E}^\c(u^\c), \lambda_\c\right>
- \eta\int\limits_{\Omega_\o} \left(Iu^\a - u^\c\right) \, dx.
\end{array}
\end{equation}
Setting the first-variations of the Lagrangian to zero yields the optimality system
\begin{equation}\label{eq:opt}
\text{find $u^\a$, $u^\c$, $\lambda_\a$, $\lambda_\c$, and $\eta$ such that $\nabla \Psi(u^\a, u^\c, \lambda_\a, \lambda_\c, \eta) = 0$,}
\end{equation}
where
\begin{equation}\label{PsiGrad}
\nabla \Psi =
\left(
 \frac{\partial \Psi}{\partial u^\a},
 \frac{\partial \Psi}{\partial u^\c},
 \frac{\partial \Psi}{\partial \lambda_\a},
\frac{\partial \Psi}{\partial \lambda_\c},
\frac{\partial \Psi}{\partial \eta}
\right)^T
\end{equation}
is the Jacobian\footnote{The notation $\frac{\partial \Psi}{\partial u^\a}$ is used to represent the vector $\frac{\partial \Psi}{\partial u^\a_\xi}$ for $\xi \in \mathcal{L}_\a$ with analogous definitions for the remaining components. } of $\Psi $.
The first-order necessary conditions (\ref{eq:opt}) are a
nonlinear system of equations for the unknowns $u^\a$, $u^\c$,
$\lambda_\a$, $\lambda_\c$, and $\eta$. To solve this system, we
employ Newton linearization. Specifically, for a given initial
guess $\mathbf{z} = [u^\a, \, \, u^\c,\, \, \lambda_\a,\,
\,\lambda_\c, \, \eta]^T$ we solve the linear equation
\begin{equation}\label{eq:newton}
\nabla^2 \Psi(\mathbf{z})\mathbf{x} = -\nabla \Psi(\mathbf{z})
\end{equation}
for the Newton increment $\mathbf{x}$ and set the new iterate to $\mathbf{z} = \mathbf{z} + \mathbf{x}$.

It is not difficult to see that the Hessian $\nabla^2 \Psi(z)$ has the form
\begin{equation}\label{PsiHess1}
\nabla^2 \Psi = \begin{pmatrix} \frac{\partial^2 \Psi}{\partial (u^\a)^2} &\frac{\partial^2 \Psi}{\partial u^\c\partial u^\a} &\frac{\partial^2 \Psi}{\partial \lambda_\a\partial u^\a} &\bm{0} &\frac{\partial^2 \Psi}{\partial \eta \partial u^\a} \\[0.3em]
\frac{\partial^2 \Psi}{\partial u^\a\partial u^\c} &\frac{\partial^2 \Psi}{\partial (u^\c)^2} &\bm{0} &\frac{\partial^2 \Psi}{\partial \lambda_\c\partial u^\c} & \frac{\partial^2 \Psi}{\partial \eta \partial u^\c}\\[0.3em]
\frac{\partial^2 \Psi}{\partial u^\a\partial \lambda_\a} &\bm{0} &\bm{0} &\bm{0}  & \bm{0} \\[0.3em]
\bm{0} &\frac{\partial^2 \Psi}{\partial u^\c\partial \lambda_\c} &\bm{0} &\bm{0}  & \bm{0} \\[0.3em]
\frac{\partial^2 \Psi}{\partial u^\a\partial \eta} &\frac{\partial^2 \Psi}{\partial u^\c\partial \eta} &\bm{0} &\bm{0} & \bm{0} \\
\end{pmatrix}
=: \begin{pmatrix} \bm{A} & \bm{B}^T \\
\bm{B} &\bm{0}
\end{pmatrix},
\end{equation}
and so, (\ref{eq:newton}) has the typical structure of a saddle-point problem.



\section{Formal Error and Complexity Analysis}\label{S:error}

We measure the error in the energy (semi-)norm,
$
\|D\bar{u}^\atc - D \bar{u}\|_{ell^2(\calL)}^2
.
$
Typically, the error has several contributions: (1) the error of truncating the infinite domain, (2) the error of modeling the atomistic interaction with the continuum interaction on a finite element mesh, and (3) an error from coupling the two models.
The first error is expected to be $\|Du\|_{\ell^2(\bbZ^d\setminus \calL)}$.
For $P_1$ (i.e., piecewise linear) elements, the second error is expected to be $\|h D^2 \bar{u} \|_{\ell^2(\calL_\c)}$, where $h$ is the element size and $D^2 u := (D_\rho D_\sigma u)_{\rho,\sigma\in\calR}$.  The third error is usually dominated by the second.
For rigorous establishments of similar error estimates, see \cite{acta.atc, OrtnerShapeev, EhrlacherOrtnerShapeev}. 
In this note, we conjecture the following result,
\begin{conj}\label{conj:error}
\begin{equation}\label{errorConj}
\|D \bar{u}^\atc - D \bar{u}\|_{\ell^2(\calL)}^2
\lesssim~
\|D \bar{u}\|_{\ell^2(\bbZ^d\setminus\calL)}^2
+
\|h D^2 \bar{u}\|_{\ell^2(\calL_\c)}^2
=:{\rm err}^2
,
\end{equation}
where $X\lesssim Y$ indicates that $X$ is less than or equal to $Y$ up to a multiplicative constant (i.e., that $\exists \, c>0$ such that $X \leq c Y$).
\end{conj}
We note that this is the most ``optimistic'' conjecture and includes only the error contributions (1) and (2) that cannot be avoided.

\subsection*{Optimal Approximation Parameters}\label{S:param}
A defect can be characterized by a far-field decay rate,
$\gamma > 0$, of the elastic displacement or stress fields.
That is, we assume that $|D^k \bar{u}(\xi)| \sim
|\xi|^{1-k-\gamma}$ (typically, $\gamma=d$ for a point defect
and $\gamma=1$ for a dislocation
\cite{mura1987micromechanics,EhrlacherOrtnerShapeev}). We
further assume a finite element discretization,
$\mathcal{T}_h$, with nodes in $\mathcal{L}_\c$ and with a radial
mesh size function $h(x) := {\rm diam}(T)$ for $x\in T$, which
will be chosen to formally optimize the error bound in~\eqref{errorConj}
subject to a fixed number of degrees of freedom.  We fully
resolve the mesh in $\Omega_\o$ so that each element of
$\mathcal{L}_\o$ is taken as a node.
The number of remaining degrees of freedom is then given by
\[
\#{\rm DoF} = \sum_{\substack{T \in \mathcal{T}_h \\ T \cap \Omega_\o = \emptyset}} 1 = \sum_{\substack{T \in \mathcal{T}_h \\ T \cap \Omega_\o = \emptyset}} \frac{|T|}{|T|} \eqsim \int_{R_\a}^{R_\c}\frac{1}{\tilde{h}^d}r^{d-1}\, {\rm d}r,
\]
where $\tilde{h}(|x|) \eqsim h(x)$ is a mesh size function that depends only on $|x|$ and $X\eqsim Y$ indicates that $X$ and $Y$ are equal up to a multiplicative constant.

Recalling that $|D^k \bar{u}(\xi)| \sim |\xi|^{1-k-\gamma}$, we
thus carry out the optimization problem:
\begin{align*}
\text{minimize}& \,
	  \int_{R_\a}^{R_\c}\tilde{h}^2 r^{-2-2\gamma}r^{d-1}\, {\rm d}r
	+ \int_{R_\c}^{\infty} r^{-2\gamma} r^{d-1}\, {\rm d}r
	\\
	\text{subject to}& \quad \left\{\begin{array}{l}\#{\rm DoF} =  \int_{R_\a}^{R_\c}\frac{1}{\tilde{h}^d}r^{d-1} \, {\rm d}r = C,
	\\ \tilde{h}(R_\a) = 1
\end{array}\right.
\end{align*}
with respect to $\tilde{h}=\tilde{h}(r)$ and $R_\c$.
Notice that we optimize only a part of the error bound, since the remaining contribution $\|h D^2 \bar{u}\|_{\ell^2(\calL_\o)} \eqsim \int_{R_\core}^{R_\a}\tilde{h}^2 r^{-2-2\gamma}r^{d-1}\, {\rm d}r$ cannot be optimized after we have fixed the mesh in $\Omega_\o$.

Introducing Lagrange multipliers and taking the variation with respect to $\tilde{h}$ we obtain
$\tilde{h}(|x|) = c |x|^{\frac{1 + \gamma}{1 + d/2}}$ for
some constant $c$, and the second constraint, $\tilde{h}(R_\a) = 1$, can then be used
to see that $h(x) = \tilde{h}(|x|) = (|x|/R_\a)^{\frac{1 +
\gamma}{1 + d/2}}$ (refer to \cite{BabuskaRheinboldt1979} for a related example of mesh optimization for ODEs).
Likewise, by differentiating with respect to $R_\c$ and using the expression for $\tilde{h}$ we find that $R_\c \eqsim R_\a^{\frac{1+\gamma}{\gamma-d/2}}$, provided $2\gamma - d > 0$.  
Finally, from the stability condition derived in~\cite{olson}, we should choose $R_\a \eqsim R_\core$.

Since the number of degrees of freedom is
\begin{equation}\label{dofEq}
{\rm DoF}
\eqsim R_\a^d +
\int_{R_\a}^{R_\c} \Big((r/R_a)^{\frac{1 +
\gamma}{1 + d/2}}\Big)^{-d} r^{d-1} {\rm d}r
\eqsim R_\a^d,
\end{equation}
we have
\begin{equation}\label{rate}
{\rm err}^2 \eqsim \left({\rm DoF}\right)^{\frac{-2-2\gamma+d}{d}}.
\end{equation}

\begin{remark}[Uniform norm]
A more involved derivation can be used to optimize the parameters for the conjecture  $\|D \bar{u}^\atc - D \bar{u}\|_{\ell^\infty(\calL)} \lesssim {\rm errinf}$, where the errors in~\eqref{errorConj} are now measured in the infinity norm.  In this case, we would get
\[
h(x) = \left(\frac{|x|}{R_\a}\right)^{1+\gamma}
,
\qquad
R_\c \eqsim R_\a^{1+\frac{1}{\gamma}}
,
\qquad
{\rm errinf} \eqsim \left({\rm DoF}\right)^{-\frac{1+\gamma}{d}}.
\]
For a dislocation (i.e., for $\gamma=1$ and $d=2$, cf.~\cite{EhrlacherOrtnerShapeev}), the energy norm is infinite so optimizing approximation parameters for $\rm err$ is ill-posed. Nevertheless, optimizing ${\rm errinf}$ is well-posed.
\end{remark}

\section{Numerical Experiments}\label{S:num}

In this section we report the results of numerical experiments conducted in 1D ($d=1$) using a next-nearest neighbor Lennard-Jones model as the underlying atomistic model.  These experiments are analogous to those run for various, popular AtC methods in~\cite{acta.atc}, with the exception that the atomistic model chosen there was the Embedded Atom Method.  Numerical experiments for the blended energy and blended force-based quasicontinuum methods using optimal approximation parameters have been presented in~\cite{bqce12,bqcf13}.  Our results provide evidence in support of the estimates conjectured in Sect.~\ref{S:error}.  We will also show how to incorporate external forces into the model as alluded to in Remark~\ref{rem:finiteDomain}.
We consider the exact, atomistic energy on the infinite lattice, $\mathbb{Z}$, to be
\begin{equation*}
\mathcal{E}^\a(u) =~ \sum_{\xi \in \mathbb{Z}} \phi(1 + D_1u(\xi)) + \phi(2 + D_1u(\xi) - D_{-1}u(\xi)) - (\phi(1) - \phi(2)) - f(\xi)u(\xi),
\end{equation*}
where $f(\xi)$ is an external force at $\xi$.  The Cauchy-Born continuum energy is
\begin{equation*}
\mathcal{E}^\c(u) =~ \int\limits_{\mathbb{R}} W(\nabla u)\, dx - \int\limits (If)u\, dx , \quad \text{where} \quad W(\mG) = \phi(1 + \mG) + \phi(2 + 2\mG),
\end{equation*}
and $If$ is the continuous linear interpolant of the force.  We assume the exact atomistic solution that we wish to approximate is (as in~\cite{acta.atc})
\[
\bar{u}^\a_{\xi} = \frac{1}{10}\left(1+\xi^2\right)^{-\gamma/2} \xi.
\]
Given this solution, we compute the external forces on an atom $\xi$ to ensure that $\bar{u}^\a$ is indeed a minimizer of the global atomistic energy.  These forces are
\[
f_{\xi} = -\frac{\partial \mathcal{E}^\a(u)}{\partial u_\xi}\big|_{u = \bar{u}^\a}.
\]
This implies the Lagrangian from Sect.~\ref{S:imp} is
\begin{equation}\label{lagrangian2}
\begin{array}{l}
\displaystyle
\Psi(u^\a, u^\c, \lambda_\a, \lambda_\c, \eta) = \frac{1}{2}\|\nabla Iu^\a - \nabla u^\c\|^{2}_{L^2(\Omega_\o)} + \left< \delta \mathcal{E}^\a(u^\a), \lambda_\a\right> \\[1ex]
\displaystyle\quad
+ \left< \delta \mathcal{E}^\c(u^\c), \lambda_\c\right>
+ \eta_1\!\! \int\limits_{\Omega_\o \cap \mathbb{R}^+} \left(Iu^\a - u^\c\right) \, dx
+ \eta_2\!\! \int\limits_{\Omega_\o \cap \mathbb{R}^-} \left(Iu^\a - u^\c\right) \, dx
\end{array}
\end{equation}
where we use the continuous, piecewise linear interpolant of $u^\a$ in this formulation.  (Recall also the need for two Lagrange multipliers to enforce the mean value zero condition on the disconnected overlap region in one dimension.)

We select a value of $R_\core$ from a range of interest, choose the mesh according to the formal analysis from Sect.~\ref{S:param}, and assign the rest of the approximation parameters via the formal derivations of Sect.~\ref{S:param}.  Namely, we set $R_\a = 2 R_\core$ and recursively construct the nodes, $\mathcal{N}_h$, of the triangulation, $\mathcal{T}_h$, as follows.  First, each $\xi \in B_{R_\a}(0)$ is chosen as a node.  Set $\xi = \max_{\zeta \in \mathcal{N}_h} \zeta$, and sequentially add a new node at $\pm\left[\xi+h(\xi)\right]$ where $h(\xi) := \lfloor (\xi/R_\a)^{\frac{1 + \gamma}{1 + d/2}}\rfloor$.  This is continued until $h(\xi) \approx \xi$, at which point we add two final nodes at $\pm R_C$.

\begin{figure}
\centering
\includegraphics[width=2.63in]{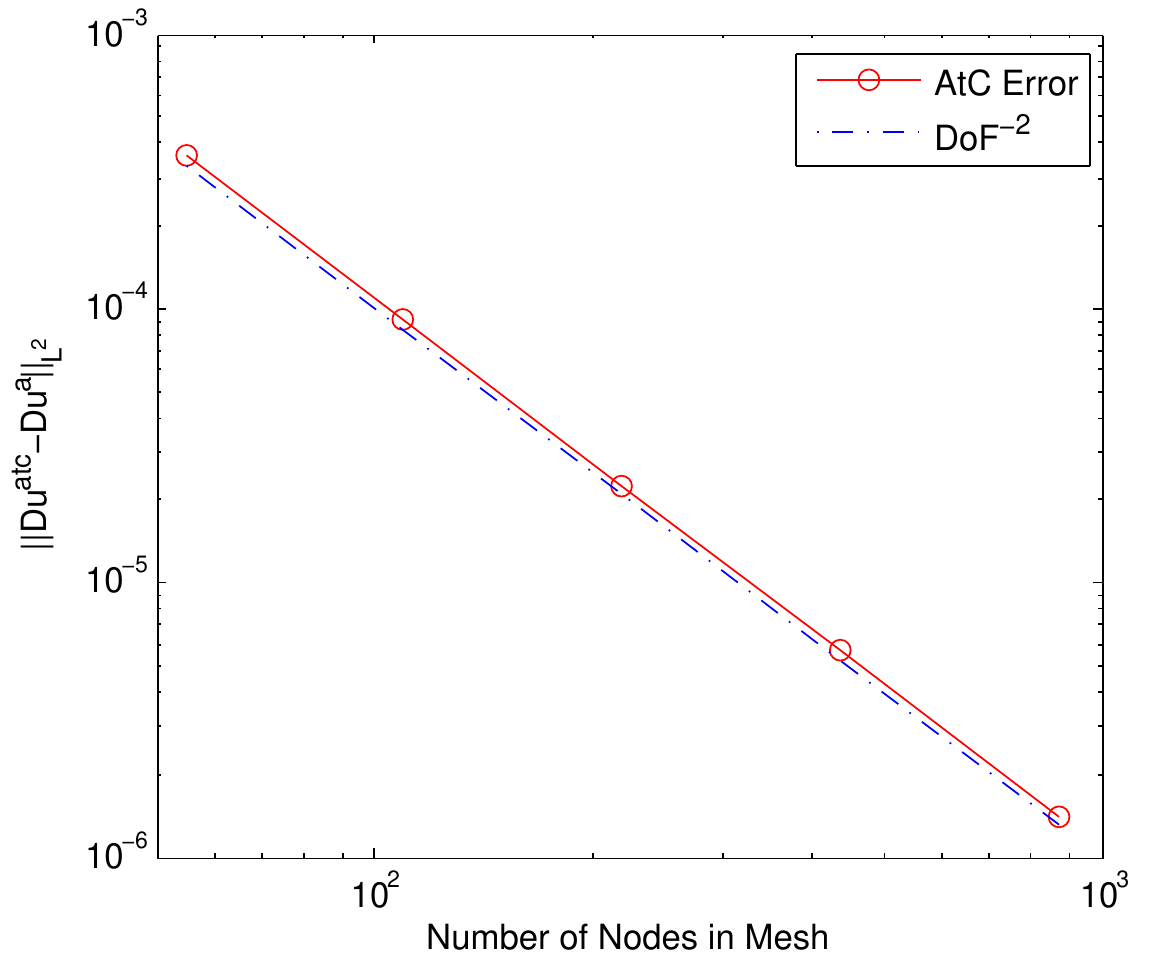}
\caption{Error of AtC approximation plotted against number of degrees of freedom.}
\label{fig:error}
\end{figure}

Finally, we take the ``defect'' approximation parameter to
be $\gamma := 3/2$ and employ our optimization-based AtC
algorithm to compute $u^\atc$ for the range of values $R_\core
\in \left\{10, \, 20, \, 40, \, 80, \, 160\right\}$.
According to our estimate~\eqref{rate} in Sect.~\ref{S:param}, we expect
the error to decay as ${\rm err}_2 \eqsim {\rm DoF}^{-2}$.  We
have plotted the error involved in each of these approximations
versus the number of degrees of freedom in Fig.~\ref{fig:error}.
In particular, the
error behaves like $(\rm DoFs)^{-2}$, which is truly optimal in
the sense of AtC methods because this is the rate of the
continuum model.  In other words, the error of coupling atomistic and continuum models is dominated by the far field error
and the continuum modeling error, as assumed in Conjecture \ref{conj:error}.

\section{Conclusion}\label{S:con}

We have formulated a new optimization-based AtC method for arbitrary interatomic potentials in multiple dimensions. Numerical simulations using a next-nearest neighbor Lennard-Jones atomistic model confirm a conjecture that the coupling error is dominated by the modeling and the domain truncation errors,
i.e., that our AtC method behaves in an optimal fashion.

\bibliographystyle{plain}	
\bibliography{NonLinOptV11}		

\end{document}